\documentclass{commat}

\title{General terms of all almost balancing numbers of first and second type}

\author{Ahmet Tekcan and Alper Erdem}

\affiliation{
    \address{Ahmet Tekcan and Alper Erdem--
    Bursa Uludag University, Faculty of Science, Department of
    Mathematics, Bursa, Turkiye
        }
    \email{tekcan@uludag.edu.tr, alper.erdem@outlook.com}
    }

\abstract{In this work, we determined the general terms of all almost ba\-lan\-cing
	numbers of first and second type in terms of balancing numbers and
	conversely we determined the general terms of all balancing numbers in terms
	of all almost balancing numbers of first and second type. We also set a
	correspondence between all almost balancing numbers of first and second type
	and Pell numbers.}



\VOLUME{31}
\YEAR{2023}
\NUMBER{1}
\firstpage{155}
\DOI{https://doi.org/10.46298/cm.10318}

\begin{paper}

\section{ Introduction}
	
	Behera and Panda (\cite{B-Pa}) defined that a positive integer $n$ is called a balancing number  if the
	Diophantine equation 
	\begin{equation}
		1+2+\cdots +(n-1)=(n+1)+(n+2)+\cdots +(n+r)  \label{bal-b11}
	\end{equation}%
	holds for some positive integer $r$ which is called balancer corresponding
	to $n$. If $n$ is a balancing number with balancer $r$, then from (\ref{bal-b11}) they get 
	\begin{equation}
r=\frac{-2n-1+\sqrt{8n^{2}+1}}{%
			2}.  \label{muty}
	\end{equation}%
	So from (\ref{muty}), they noted that $n$ is a balancing number if and only if $8n^{2}+1$ is a perfect square. Though the
	definition of balancing num\-bers suggests that no balancing number should
	be less than $2$. But from (\ref{muty}), they noted that $8(0)^{2}+1=1$ and $%
	8(1)^{2}+1=3^{2}$ are perfect squares. So they accepted that $0$
	and $1$ to be balancing numbers. Let $B_{n}$ denote the $n^{\text{th}}$
	balancing number. Then $B_{0}=0,B_{1}=1,B_{2}=6$ and $B_{n+1}=6B_{n}-B_{n-1}$
	for $n\geq 2$.
	
	Later Panda and Ray (\cite{pa-ray}) defined that a positive integer $n$ is
	called a cobalancing number if the Diophantine equation 
	\begin{equation}
		1+2+\cdots +n=(n+1)+(n+2)+\cdots +(n+r)  \label{bal-b1}
	\end{equation}%
	holds for some positive integer $r$ which is called cobalancer corresponding
	to $n$. If $n$ is a cobalancing number with cobalancer $r$, then from (\ref{bal-b1}) they get 
	\begin{equation}
r=\frac{-2n-1+\sqrt{%
				8n^{2}+8n+1}}{2}.  \label{muty1}
	\end{equation}%
	So from (\ref{muty1}), they noted that $n$ is a cobalancing number if and only if $8n^{2}+8n+1$ is a perfect square.
	Since $8(0)^{2}+8(0)+1=1$ is a perfect square, they accepted $0$ to be a
	cobalancing number just like Behera and Panda accepted $0$ and $1$ to be
	balancing numbers. Let $b_{n}$ denote the $n^{\text{th}}$ cobalancing number.
	Then $b_{0}=b_{1}=0,b_{2}=2$ and $b_{n+1}=6b_{n}-b_{n-1}+2$ for $n\geq 2$.
	
	It is clear from (\ref{bal-b11}) and (\ref{bal-b1}) that every ba\-lancing
	number is a cobalancer and every cobalancing number is a balancer, that is, $%
	B_{n}=r_{n+1}$ and $R_{n}=b_{n}$ for $n\geq 1$, where $R_{n}$ is the $n^{%
		\text{th}}$ the balancer and $r_{n}$ is the $n^{\text{th}}$ cobalancer.
	Since $R_{n}=b_{n}$, we get from (\ref{bal-b11}) that 
	\begin{equation}
		b_{n}=\frac{-2B_{n}-1+\sqrt{8B_{n}^{2}+1}}{2}\text{ \ and\ \ }B_{n}=\frac{b_{n}+1+\sqrt{8b_{n}^{2}+8b_{n}+1}}{2}.\label{baa-12}
	\end{equation}%
	Thus from (\ref{baa-12}), $B_{n}$ is a balancing number if and
	only if $8B_{n}^{2}+1$ is a perfect square and $b_{n}$ is a cobalancing
	number if and only if $8b_{n}^{2}+8b_{n}+1$ is a perfect square. So 
	\begin{equation}
		C_{n}=\sqrt{8B_{n}^{2}+1}\ \ \text{and}\ \ c_{n}=\sqrt{8b_{n}^{2}+8b_{n}+1}
		\label{lucaslar}
	\end{equation}%
	are integers which are called the $n^{\text{th}}$ Lucas-balancing number and 
	$n^{\text{th}}$ Lucas-cobalancing number, respectively.
	
	Let $\alpha =1+\sqrt{2}$ and $\beta =1-\sqrt{2}$ be the roots of the
	characteristic equation for Pell numbers $P_{n}$. Then Binet
	formulas for balancing numbers, co\-ba\-lan\-cing numbers, Lucas-balancing numbers
	and Lucas-cobalancing numbers are $$B_{n}=\frac{\alpha ^{2n}-\beta ^{2n}}{4%
		\sqrt{2}},b_{n}=\frac{\alpha ^{2n-1}-\beta ^{2n-1}}{4\sqrt{2}}-\frac{1}{2}%
	,C_{n}=\frac{\alpha ^{2n}+\beta ^{2n}}{2}\ \ \text{and} \ \ c_{n}=\frac{\alpha
		^{2n-1}+\beta ^{2n-1}}{2}$$ for $n\geq 1$, respectively (see also \cite	{tkcn1}, \cite{olajas1}, \cite{pa-ray12}, \cite{raytez}, \cite{tkcn2}).
	
	Balancing numbers and their generalizations have been investigated by
	se\-ve\-ral authors from many aspects. In \cite{lip}, Liptai proved that
	there is no Fi\-bo\-nac\-ci ba\-lan\-cing number except $1$ and in \cite%
	{lip1} he proved that there is no Lucas balancing number. In \cite{szalay},
	Szalay considered the same problem and obtained some nice results by a
	different method. In \cite{tunde}, Kov\'{a}cs, Liptai and Olajos
	ex\-ten\-ded the concept of balancing numbers to the $(a,b)$-balancing
	numbers defined as follows: Let $a>0$ and $b\geq 0$ be coprime integers. If 
	\begin{equation*}
		(a+b)+\cdots +(a(n-1)+b)=(a(n+1)+b)+\cdots +(a(n+r)+b)
	\end{equation*}%
	for some positive integers $n$ and $r$, then $an+b$ is an $(a,b)$-balancing
	number. The sequence of $(a,b)$-balancing numbers is denoted by $%
	B_{m}^{(a,b)}$ for $m\geq 1$. In \cite{akos}, Liptai, Luca, Pint\'{e}r and
	Szalay generalized the notion of balancing numbers to numbers defined as
	follows: Let $y,k,l\in \mathbb{Z}^{+}$ such that $y\geq 4$. Then a po\-si\-ti\-ve
	integer $x$ with $x\leq y-2$ is called a $(k,l)$-power numerical center for $%
	y$ if $$1^{k}+\cdots +(x-1)^{k}=(x+1)^{l}+\cdots +(y-1)^{l}.$$ They studied
	the number of solutions of the equation above and proved several effective
	and ineffective finiteness results for $(k,l)$-power numerical centers. For
	positive integers $k,x$, let $$\Pi _{k}(x)=x(x+1)\dots (x+k-1).$$ Then it was
	proved in \cite{tunde} that the equation $B_{m}=\Pi _{k}(x)$ for fixed
	integer $k\geq 2$ has only infinitely many solutions and for $k\in \{2,3,4\}$
	all solutions were determined. In \cite{tengely} Tengely con\-si\-de\-red
	the case $$B_{m}=x(x+1)(x+2)(x+3)(x+4)$$ for $k=5$ and proved that this
	Diophantine equation has no solution for $m\geq 0$ and $x\in \mathbb{Z}$. In 
	\cite{komats}, Panda, Ko\-mat\-su and Davala considered the reciprocal sums
	of sequences involving balancing and Lucas-balancing numbers. In \cite{patel1}, Patel, Irmak and Ray considered incomplete balancing and Lucas-balancing
	numbers and in \cite{raysums}, Ray considered the sums of ba\-lan\-cing and
	Lucas-balancing numbers by matrix methods. In \cite{erdem}, Tekcan and Erdem  considered the $t$-cobalancing numbers and $t$-cobalancers, in \cite{ayd}, Tekcan and Ayd\i n considered the $t$-balancers, $t$-balancing numbers and Lucas $t$-balancing numbers and in \cite{meryem}, Tekcan and Y\i ld\i z considered the balcobalancing numbers and balcobalancers.

	\section{Results}
	
	In this section we determine the general terms of almost balancing
	numbers, almost cobalancing
	numbers, almost Lucas-balancing
	numbers and almost Lucas-cobalancing
	numbers of first and second type. Almost balancing numbers first
	defined by Panda and Panda in \cite{almost}. A positive integer $n$ is
	called an almost ba\-lan\-cing num\-ber if the Di\-op\-han\-ti\-ne equation 
	\begin{equation}
		\left\vert \lbrack (n+1)+(n+2)+\cdots +(n+r)]-[1+2+\cdots +(n-1)]\right\vert
		=1 \label{muyer}
	\end{equation}%
	holds for some positive integer $r$ which is called the almost balancer. 
	
	From (\ref{muyer}), they have two cases: If $[(n+1)+(n+2)+\cdots +(n+r)]-[1+2+\cdots +(n-1)]=1$, then $n$ is called an almost balancing
	number of first type and $r$ is called an almost balancer of first type and
	in this case 
	\begin{equation}
		r=\frac{-2n-1+\sqrt{8n^{2}+9}}{2}.  \label{alm1-23}
	\end{equation}%
	If $[(n+1)+(n+2)+\cdots +(n+r)]-[1+2+\cdots +(n-1)]=-1$, then $n$ is called an almost
	balancing number of second type and $r$ is called an almost balancer of
	second type and in this case 
	\begin{equation}
		r=\frac{-2n-1+\sqrt{8n^{2}-7}}{2}.  \label{alm1-13}
	\end{equation}%

	Let $B_{n}^{\ast }$ and $B_{n}^{\ast \ast }$ denote the $n^{\text{th}}$
	almost balancing number of first type and of second type, respectively. Then from (\ref{alm1-23}), $%
	B_{n}^{\ast }$ is an almost balancing number of first type if and only if $%
	8(B_{n}^{\ast })^{2}+9$ is a perfect square and from (\ref{alm1-13}), $%
	B_{n}^{\ast \ast }$ is an almost balancing number of second type if and only
	if $8(B_{n}^{\ast \ast })^{2}-7$ is a perfect square. Thus 
	\begin{equation}
		C_{n}^{\ast }=\sqrt{8(B_{n}^{\ast })^{2}+9}\text{ \ and\ \ }C_{n}^{\ast \ast
		}=\sqrt{8(B_{n}^{\ast \ast })^{2}-7}  \label{alm1-2}
	\end{equation}%
	are integers which are called the $n^{\text{th}}$ almost Lucas-balancing
	number of first type and of second type, respectively.
	
	Later in \cite{almosttez}, Panda defined that a positive integer $n$ is
	called an almost cobalancing number if the Di\-op\-han\-ti\-ne equation 
	\begin{equation}
		\left\vert \lbrack (n+1)+(n+2)+\cdots +(n+r)]-(1+2+\cdots +n)\right\vert =1\label{muty12}
	\end{equation}%
	holds for some positive integer $r$ which is called an almost cobalancer. 
	
	From (\ref{muty12}), they have two cases: If $[(n+1)+(n+2)+\cdots +(n+r)]-(1+2+\cdots +n)=1$, then $n$ is called an almost
	cobalancing number of first type and $r$ is called an almost cobalancer of
	first type and in this case 
	\begin{equation}
		r=\frac{-2n-1+\sqrt{8n^{2}+8n+9}}{2}.  \label{asut12}
	\end{equation}%
	If $[(n+1)+(n+2)+\cdots +(n+r)]-(1+2+\cdots +n)=-1$, then $n$ is called an almost
	cobalancing number of second type and $r$ is called an almost cobalancer of
	second type and in this case 
	\begin{equation}
		r=\frac{-2n-1+\sqrt{8n^{2}+8n-7}}{2}.  \label{asut1}
	\end{equation}%

	Let $b_{n}^{\ast }$ and $b_{n}^{\ast \ast }$ denote the $n^{\text{th}}$
	almost cobalancing number of first type and of second type, respectively. Then from (\ref{asut12}), $%
	b_{n}^{\ast }$ is an almost cobalancing number of first type if and only if $%
	8(b_{n}^{\ast })^{2}+8b_{n}^{\ast }+9$ is a perfect square and from (\ref%
	{asut1}), $b_{n}^{\ast \ast }$ is an almost cobalancing number of second
	type if and only if $8(b_{n}^{\ast \ast })^{2}+8b_{n}^{\ast \ast }-7$ is a
	perfect square. Thus%
	\begin{equation}
		c_{n}^{\ast }=\sqrt{8(b_{n}^{\ast })^{2}+8b_{n}^{\ast }+9}\text{ \ and\ \ }%
		c_{n}^{\ast \ast }=\sqrt{8(b_{n}^{\ast \ast })^{2}+8b_{n}^{\ast \ast }-7}
		\label{alm1-3}
	\end{equation}%
	are integers which are called the $n^{\text{th}}$ almost Lucas-cobalancing
	number of first type and of second type, res\-pec\-ti\-vely.

	Like in balancing numbers, we notice that every almost
	ba\-lan\-cing number is an almost cobalancer and every almost cobalancing
	number is an almost balancer, that is, $B_{n}^{\ast }=r_{n+1}^{\ast \ast
	},$ $B_{n}^{\ast \ast }=r_{n}^{\ast },b_{n}^{\ast }=R_{n+2}^{\ast }$ and $%
	b_{n}^{\ast \ast }=R_{n}^{\ast }$ for $n\geq 1$, where $R_{n}^{\ast }$ is
	the $n^{\text{th}}$ almost balancer of first type, $R_{n}^{\ast \ast }$ is
	the $n^{\text{th}}$ almost balancer of second type, $r_{n}^{\ast }$ is the $%
	n^{\text{th}}$ almost cobalancer of first type and $r_{n}^{\ast \ast }$ is
	the $n^{\text{th}}$ almost cobalancer of second type.
	
	\subsection{Almost Balancing and Almost Lucas-Balancing Numbers of First and Second Type.}
	
	We see in (\ref%
	{alm1-23}) that $x$ is an almost balancing number of first type if and only
	if $8x^{2}+9$ is a perfect square and in (\ref{alm1-13}), $x$ is an almost
	ba\-lancing number of second type if and only if $8x^{2}-7$ is a perfect
	square. Let $8x^{2}+9=y^{2}$ and let $8x^{2}-7=w^{2}$ for some positive
	integers $y$ and $w.$ Then we get the Pell equations (\cite	{barb}, \cite{Jacob}) 
	\begin{equation}
		8x^{2}-y^{2}=-9\text{ \ and \ }8x^{2}-w^{2}=7.  \label{pell1}
	\end{equation}%
	For the set of all (positive) integer solutions of (\ref{pell1}), we need some notations: Let $\Delta $ be a
	non-square discriminant. Then the $\Delta $-order $O_{\Delta }$ is defined
	to be the ring
	\[
	O_{\Delta }=\{x+y\rho _{\Delta }:x,y\in \mathbb{Z\}},
	\]
	 where 
	$\rho _{\Delta }=\sqrt{\frac{\Delta }{4}}$ if $\Delta \equiv 0(\text{mod}\
	4) $ or $\frac{1+\sqrt{\Delta }}{2}$ if $\Delta \equiv 1(\text{mod}\ 4)$. So 
	$O_{\Delta }$ is a subring of $\mathbb{Q(}\sqrt{\Delta })\mathbb{=}\{x+y%
	\sqrt{\Delta }:x,y\in \mathbb{Q}\}$. The unit group $O_{\Delta }^{u}$ is
	defined to be the group of units of the ring $O_{\Delta }^{{}}$. We can rewrite an integral indefinite quadratic form (\cite{flath}) $F(x,y)=ax^{2}+bxy+cy^{2}$ of discriminant $\Delta$ to be
	$$F(x,y)=\frac{(xa+y\frac{b+\sqrt{\Delta }}{2})(xa+y\frac{b-\sqrt{\Delta }}{2})%
	}{a}.$$ So the module $M_{F}$ of $F$ is $M_{F}=\{xa+y\frac{b+\sqrt{\Delta }}{2}:x,y\in \mathbb{Z}\}\subset \mathbb{Q}(%
	\sqrt{\Delta }).$ Therefore we get $(u+v\rho _{\Delta })(xa+y\frac{b+\sqrt{\Delta }}{2}%
	)=x^{\prime }a+y^{\prime }\frac{b+\sqrt{\Delta }}{2}$, where 
	\begin{equation}
		\lbrack x^{\prime }\ \ y^{\prime }]=\left\{ 
		\begin{array}{cc}
			\lbrack x\ \ y]\left[ 
			\begin{array}{cc}
				u-\frac{b}{2}v & av \\ 
				-cv & u+\frac{b}{2}v%
			\end{array}%
			\right] \text{ \ \ } & \text{if }\Delta \equiv 0(\text{mod}\ 4) \\ 
			\lbrack x\ \ y]\left[ 
			\begin{array}{cc}
				u+\frac{1-b}{2}v & av \\ 
				-cv & u+\frac{1+b}{2}v%
			\end{array}%
			\right] & \text{if }\Delta \equiv 1(\text{mod}\ 4).%
		\end{array}%
		\right.  \label{mat1x}
	\end{equation}%
	Let $m$ be any integer and let $\Omega $ denote the set of all integer
	solutions of $F(x,y)=m$, that is, $\Omega =\{(x,y):F(x,y)=m\}$. Then there
	is a bijection $\Psi :\Omega \rightarrow \{\gamma \in M_{F}:N(\gamma )=am\}.$
	The action of $O_{\Delta ,1}^{u}=\{\alpha \in O_{\Delta }^{u}:N(\alpha )=1\}$
	on the set $\Omega $ is most interesting when $\Delta $ is a positive
	non-square since $O_{\Delta ,1}^{u}$ is infinite. Therefore the orbit of
	each solution will be infinite and so the set $\Omega $ is either empty or
	infinite. Since $O_{\Delta ,1}^{u}$ can be explicitly determined, the set $%
	\Omega $ is satisfactorily described by the representation of such a list,
	called a set of representatives of the orbits. Let $\varepsilon _{\Delta }$
	be the smallest unit of $O_{\Delta }^{{}}$ that is grater than $1$ and let $%
	\tau _{\Delta }=\varepsilon _{\Delta }$ if $N(\varepsilon _{\Delta })=1$ or $%
	\varepsilon _{\Delta }^{2}$ if $N(\varepsilon _{\Delta })=-1$. Then every $%
	O_{\Delta ,1}^{u}$ orbit of integral solutions of $F(x,y)=m$ contains a
	solution $(x,y)\in \mathbb{Z}\times \mathbb{Z}$ such that $0\leq y\leq U$,
	where $U=\left\vert \frac{am\tau _{\Delta }}{\Delta }\right\vert ^{\frac{1}{2%
	}}(1-\frac{1}{\tau _{\Delta }})$ if $am>0$ or $U=\left\vert \frac{am\tau
		_{\Delta }}{\Delta }\right\vert ^{\frac{1}{2}}(1+\frac{1}{\tau _{\Delta }})$
	if $am<0$. So for finding a set of representatives of the $O_{\Delta ,1}^{u}$
	orbits of integral solutions of $F(x,y)=m$, we must find for each integer $%
	y_{0}$ in the range $0\leq y_{0}\leq U$, whether $\Delta y_{0}^{2}+4am$ is a
	perfect square or not since $$ax_{0}^{2}+bx_{0}y_{0}+cy_{0}^{2}=m\Leftrightarrow \Delta y_{0}^{2}+4am=(2ax_{0}+by_{0})^{2}.$$ If $\Delta
	y_{0}^{2}+4am$ is a perfect square, then $x_{0}=\frac{-by_{0}\pm \sqrt{%
			\Delta y_{0}^{2}+4am}}{2a}.$ So there is a set of representatives $\text{Rep}%
	=\{[x_{0}\ \ y_{0}]\}$. Thus for the matrix $M$ defined in (\ref{mat1x}),
	the set of all integer solutions of $F(x,y)=m$ is $\Omega=\{\pm (x,y):[x\ \
	y]=[x_{0}\ \ y_{0}]M^{n},n\in \mathbb{Z}\}$.

	For the set of all integer solutions of (\ref{pell1}), we can can give the following
	theorem.
	
	\begin{theorem}
		\label{teo1} The set of all integer solutions of $8x^{2}-y^{2}=-9$ is $\Omega=\{(3B_{n},3C_{n}):n\geq 1\},$ and the set of all integer solutions of $8x^{2}-w^{2}=7$ is 
		\[
		\Omega=\{(B_{n-1}+C_{n-1},8B_{n-1}+C_{n-1}):n\geq 1\}\cup \{(-B_{n}+C_{n},8B_{n}-C_{n}):n\geq 1\}.
		\]
	\end{theorem}
	
	\begin{proof}For the Pell equation $8x^{2}-y^{2}=-9$, we have $F(x,y)=8x^{2}-y^{2}$
		of discriminant $\Delta =32$. So we get $\tau _{32}=3+\sqrt{8}$. Thus the
		set of representatives is $\text{Rep}=\{[0\ \ 3]\}$ and $M=\left[ 
		\begin{array}{cc}
			3 & 8 \\ 
			1 & 3%
		\end{array}%
		\right] $ by (\ref{mat1x}). Here we notice that $[0\ \ 3]M^{n}$ generates all integer
		solutions $(x_{n},y_{n})$ of $8x^{2}-y^{2}=-9$ for $n\geq 1$. It can be easily seen that the $n^{\text{th}}$ power of $M$ is 
		\begin{equation*}
			M^{n}=\left[ 
			\begin{array}{cc}
				C_{n} & 8B_{n} \\ 
				B_{n} & C_{n}%
			\end{array}%
			\right] 
		\end{equation*}%
		for $n\geq 1$. Thus the set of all integer solutions is $\Omega=\{(3B_{n},3C_{n}):n\geq 1\}$.
		
		For the second Pell equation $8x^{2}-w^{2}=7$, we get $\tau _{32}=3+\sqrt{8}$. So the set of representatives is $\text{Rep}=\{[\pm 1\ \ 1]\}$ and in this case $%
		[1\ \ \ 1]M^{n-1}$ generates all in\-te\-ger solutions $(x_{2n-1},w_{2n-1})$ and 
		$[1\ \ -1]M^{n}$ generates all integer solutions $(x_{2n},w_{2n})$ for $%
		n\geq 1$. Thus the result is obvious. 
	\end{proof}
	
	From Theorem \ref{teo1}, we can give the following theorem. 
	
	\begin{theorem}
		\label{teo2}The general terms of almost balancing and almost Lucas-balancing
		numbers of first type are 
		\begin{equation*}
			B_{n}^{\ast }=3B_{n},\ C_{n}^{\ast }=3C_{n}
		\end{equation*}%
		for $n\geq 1,$ and the general terms of almost balancing and almost
		Lucas-ba\-lan\-cing numbers of second type are%
		\begin{eqnarray*}
			&&B_{2n-1}^{\ast \ast }=B_{n-1}+C_{n-1},\ B_{2n}^{\ast \ast }=-B_{n}+C_{n} \\
			&&C_{2n-1}^{\ast \ast }=8B_{n-1}+C_{n-1},\ C_{2n}^{\ast \ast }=8B_{n}-C_{n}
		\end{eqnarray*}%
		for $n\geq 1.$
	\end{theorem}
	
	\begin{proof}
		We proved in Theorem \ref{teo1} that the set of all integer
		solutions of $8x^{2}-y^{2}=-9$ is $\Omega=\{(3B_{n},3C_{n}):n\geq 1\}$. Since $%
		x=B_{n}^{\ast }$, we get $B_{n}^{\ast }=3B_{n}$. So from
		(\ref{alm1-2}), we deduce that 
		\begin{equation*}
			C_{n}^{\ast }=\sqrt{8(B_{n}^{\ast })^{2}+9}=\sqrt{8(3B_{n})^{2}+9}=3\sqrt{8B_{n}^{2}+1}=3C_{n}.
		\end{equation*}%
		Similarly since the set of all integer
		solutions of $8x^{2}-w^{2}=7$ is
		\[
		\Omega=\{(B_{n-1}+C_{n-1},8B_{n-1}+C_{n-1}):n\geq 1\}\cup \{(-B_{n}+C_{n},8B_{n}-C_{n}):n\geq 1\},
		\]
		we get $B_{2n-1}^{\ast
			\ast}=B_{n-1}+C_{n-1},$ $C_{2n-1}^{\ast \ast
		}=8B_{n-1}+C_{n-1},B_{2n}^{\ast \ast }=-B_{n}+C_{n}$ and
		$C_{2n}^{\ast \ast }=8B_{n}-C_{n}$ for $n\geq 1$.
	\end{proof}

	Here we note that $B_{0}^{\ast }=0,C_{0}^{\ast }=3,B_{0}^{\ast \ast }=1,C_{0}^{\ast \ast }=-1$. Also since $8(1)^{2}-7=1$ and $8(2)^{2}-7=5^{2}$ are
	perfect squares by (\ref{alm1-13}), we accept $1$ and $2$ be almost balancing numbers of second type.

	\subsection{Almost Cobalancing and Almost Lucas-Cobalancing Numbers of First and Second Type.}
	
	In this subsection, we will determine the general terms of almost
	co\-ba\-lan\-cing and almost Lucas-co\-ba\-lan\-cing numbers of first and
	second type. Since $n$ is an almost cobalancing number of first type if and
	only if $8n^{2}+8n+9$ is a perfect square by (\ref{asut12}) and $n$ is an
	almost cobalancing number of second type if and only if $8n^{2}+8n-7$ is a
	perfect square by (\ref{asut1}), we set $8n^{2}+8n+9=y^{2}$ and $%
	8n^{2}+8n-7=w^{2}$ for some positive integers $y$ and $w$. Then we get the
	equations $2(2n+1)^{2}-y^{2}=-7$ and $2(2n+1)^{2}-w^{2}=9$. Taking $2n+1=x$,
	we get the Pell equations 
	\begin{equation}
		2x^{2}-y^{2}=-7\text{ \ and \ }2x^{2}-w^{2}=9.  \label{pell2}
	\end{equation}%
	For the set of all integer solutions of (\ref{pell2}), we can can give the following
	theorem.
	
	\begin{theorem}
		\label{teo3} The set of all integer
		solutions of $2x^{2}-y^{2}=-7$ is
		\[
		\Omega=\{(6B_{n-1}+C_{n-1}, 4B_{n-1}+3C_{n-1}):n\geq 1\}\cup \{(6B_{n}-C_{n},-4B_{n}+3C_{n}):n\geq 1\},
		\]
		and the set of all integer
		solutions of $2x^{2}-w^{2}=9$ is 
		\[
		\Omega=\{(6B_{n-1}+3C_{n-1},12B_{n-1}+3C_{n-1}):n\geq 1\}.
		\]
	\end{theorem}
	
	\begin{proof}For the Pell equation $2x^{2}-y^{2}=-7$, we get $F(x,y)=2x^{2}-y^{2}$ of
		discriminant $\Delta =8$. So $\tau _{8}=3+2\sqrt{2}$ and hence the set
		of representatives is $\text{Rep}=\{[\pm 1\ \ 3]\}$ and $M=\left[ 
		\begin{array}{cc}
			3 & 4 \\ 
			2 & 3%
		\end{array}%
		\right] $. Here $[1\ \ 3]M^{n-1}$ ge\-ne\-ra\-tes all integer solutions $%
		(x_{2n-1},y_{2n-1})$ and $[-1\ \ 3]M^{n}$ generates all integer solutions $%
		(x_{2n},y_{2n})$ for $n\geq 1$. Since the $n^{\text{th}}$ power of $M$ is 
		\begin{equation*}
			M^{n}=\left[ 
			\begin{array}{cc}
				C_{n} & 4B_{n} \\ 
				2B_{n} & C_{n}%
			\end{array}%
			\right]
		\end{equation*}%
		for $n\geq 1$, we deduce that the set of all integer
		solutions is 
		\[
		\Omega=\{(6B_{n-1}+C_{n-1},4B_{n-1}+3C_{n-1}):n\geq 1\}\cup \{(6B_{n}-C_{n},-4B_{n}+3C_{n}):n\geq		1\}.
		\]
		
		For the second Pell equation $2x^{2}-w^{2}=9$, we get $\tau _{8}=3+2\sqrt{2%
		}$ and the set of representatives is $\text{Rep}=\{[\pm 3\ \ 3]\}$. In this case $[3\ \
		3]M^{n-1}$ generates all integer so\-lu\-ti\-ons $(x_{n},w_{n})$ for $n\geq 1
		.$ Thus the result is obvious. 
	\end{proof}
	
	From Theorem \ref{teo3}, we can give the following theorem.

	\begin{theorem}
		\label{teo4}The general terms of almost cobalancing and almost
		Lucas-co\-ba\-lan\-cing numbers of first type are 
		\begin{eqnarray*}
			&&b_{2n}^{\ast }=2b_{n+1}-b_{n},\ b_{2n-1}^{\ast }=4b_{n}-b_{n-1}+1, \\
			&&c_{2n}^{\ast }=c_{n+2}-4c_{n+1},\ c_{2n-1}^{\ast }=c_{n+1}-2c_{n}
		\end{eqnarray*}%
		for $n\geq 1,$ and the general terms of almost cobalancing and almost
		Lucas-co\-ba\-lan\-cing numbers of second type are%
		\begin{equation*}
			b_{n}^{\ast \ast }=3b_{n}+1,\ c_{n}^{\ast \ast }=3c_{n}
		\end{equation*}%
		for $n\geq 1.$
	\end{theorem}
	
	\begin{proof}
		We proved in Theorem \ref{teo3} that the set of all integer
		solutions of $2x^{2}-y^{2}=-7$ is $\Omega=\{(6B_{n-1}+C_{n-1},4B_{n-1}+3C_{n-1}):n\geq 1\}\cup \{(6B_{n}-C_{n},-4B_{n}+3C_{n}):n\geq
		1\}$. Since $x=2n+1$, we get
		\begin{align*}
			b_{2n}^{\ast }& =\frac{6B_{n}+C_{n}-1}{2} \\
			& =\frac{6(\frac{\alpha ^{2n}-\beta
					^{2n}}{4\sqrt{2}})+\frac{\alpha
					^{2n}+\beta ^{2n}}{2}-1}{2}\\
			& =\frac{\alpha ^{2n}(2\alpha -\alpha ^{-1})+\beta ^{2n}(-2\beta
				+\beta ^{-1})}{4\sqrt{2}}-\frac{1}{2}\\
			& =2(\frac{\alpha ^{2n+1}-\beta ^{2n+1}}{4\sqrt{2}}-\frac{1}{2})-(\frac{%
				\alpha ^{2n-1}-\beta ^{2n-1}}{4\sqrt{2}}-\frac{1}{2}) \\
			& =2b_{n+1}-b_{n}.
		\end{align*}%
		Thus from (\ref{alm1-3}), we get
		\begin{align*}
			c_{2n}^{\ast }& =\sqrt{8(b_{2n}^{\ast })^{2}+8b_{2n}^{\ast }+9} \\
			& =\sqrt{8(2b_{n+1}-b_{n})^{2}+8(2b_{n+1}-b_{n})+9} \\
			& =\sqrt{\alpha ^{4n}(\frac{11+6\sqrt{2}}{4})+\beta ^{4n}(\frac{11-6\sqrt{2}%
				}{4})+\frac{7}{2}} \\
			& =\sqrt{%
				\begin{array}{c}
					(\frac{\alpha ^{2n+3}+\beta ^{2n+3}}{2})^{2}-4(\frac{\alpha
						^{2n+3}+\beta ^{2n+3}}{2})(\alpha ^{2n+1}+\beta ^{2n+1})+4(\alpha ^{2n+1}+\beta ^{2n+1})^{2}%
				\end{array}%
			} \\
			& =\sqrt{[(\frac{\alpha ^{2n+3}+\beta ^{2n+3}}{2}-4(\frac{\alpha
					^{2n+1}+\beta ^{2n+1}}{2})]^{2}} \\
			& =\frac{\alpha ^{2n+3}+\beta ^{2n+3}}{2}-4(\frac{\alpha
				^{2n+1}+\beta
				^{2n+1}}{2}) \\
			& =c_{n+2}-4c_{n+1}.
		\end{align*}%
		The others can be proved similarly.
	\end{proof}
	
	Here we note that $b_{0}^{\ast }=0,c_{0}^{\ast }=3,$ $b_{0}^{\ast \ast
	}=1$ and $c_{0}^{\ast \ast }=3$. Also since $8(1)^{2}+8(1)$ $-7=3^{2}$ is a perfect
	square by (\ref{asut1}), we accept $1$ to be an almost cobalancing number of second type.

	In Theorems \ref{teo2} and \ref{teo4}, we deduce the general terms of all
	almost ba\-lan\-cing numbers of first and second type in terms of balancing
	numbers. Conversely, we can deduce the general terms of all balancing
	numbers in terms of all almost balancing numbers of first and second type as
	follows:
	
	\begin{theorem}
		\label{teo5}The general terms of all balancing numbers are 
		\begin{equation*}
			B_{n}=\frac{B_{n}^{\ast }}{3},\ b_{n}=\frac{b_{2n-1}^{\ast }-b_{2n-2}^{\ast
				}-1}{2},\ C_{n}=\frac{C_{n}^{\ast }}{3},\ c_{n}=\frac{c_{2n-1}^{\ast
				}-c_{2n-2}^{\ast }}{2}
		\end{equation*}%
		for $n\geq 1,$ or 
		\begin{equation*}
			B_{n}=\frac{B_{2n+1}^{\ast \ast }-B_{2n}^{\ast \ast }}{2},\ b_{n}=\frac{%
				b_{n}^{\ast \ast }-1}{3},\ C_{n}=\frac{C_{2n+1}^{\ast \ast }-C_{2n}^{\ast
					\ast }}{2},\ c_{n}=\frac{c_{n}^{\ast \ast }}{3}
		\end{equation*}%
		for $n\geq 1.$
	\end{theorem}
	
	\begin{proof}
		The result is obvious from Theorems \ref{teo2} and \ref{teo4}.
	\end{proof}
	
	Thus we construct a one-to-one correspondence between all balancing numbers
	and all almost balancing numbers of first and second type. Moreover, the general terms
	of all almost balancing numbers of  first type can be given in terms of all almost balancing numbers
	of second type and conversely the general terms of all almost balancing numbers of
	second type can be given in terms of all almost balancing numbers of first type as follows.

	\begin{theorem}
		\label{teo6}The general terms of all almost balancing numbers of first type
		are 
		\begin{eqnarray*}
			&&B_{n}^{\ast }=\frac{3B_{2n+1}^{\ast \ast }-3B_{2n}^{\ast \ast }}{2},\
			C_{n}^{\ast }=\frac{3C_{2n+1}^{\ast \ast }-3C_{2n}^{\ast \ast }}{2}, \\
			&&b_{2n-1}^{\ast }=\frac{4b_{n}^{\ast \ast }-b_{n-1}^{\ast \ast }}{3},\
			b_{2n}^{\ast }=\frac{2b_{n+1}^{\ast \ast }-b_{n}^{\ast \ast }-1}{3}, \\
			&&c_{2n-1}^{\ast }=\frac{c_{n+1}^{\ast \ast }-2c_{n}^{\ast \ast }}{3},\
			c_{2n}^{\ast }=\frac{c_{n+2}^{\ast \ast }-4c_{n+1}^{\ast \ast }}{3}
		\end{eqnarray*}%
		for $n\geq 1,$ and the general terms of all almost balancing numbers of
		second type are 
		\begin{eqnarray*}
			&&B_{2n-1}^{\ast \ast }=\frac{B_{n-1}^{\ast }+C_{n-1}^{\ast }}{3},\
			B_{2n}^{\ast \ast }=\frac{-B_{n}^{\ast }+C_{n}^{\ast }}{3}, \\
			&&b_{n}^{\ast \ast }=\frac{3b_{2n-1}^{\ast }-3b_{2n-2}^{\ast }-1}{2},\
			c_{n}^{\ast \ast }=\frac{3c_{2n-1}^{\ast }-3c_{2n-2}^{\ast }}{2}, \\
			&&C_{2n-1}^{\ast \ast }=\frac{8B_{n-1}^{\ast }+C_{n-1}^{\ast }}{3},\
			C_{2n}^{\ast \ast }=\frac{8B_{n}^{\ast }-C_{n}^{\ast }}{3}
		\end{eqnarray*}%
		for $n\geq 1$.
	\end{theorem}
	
	\begin{proof}
		Since $B_{n}^{\ast }=3B_{n}$ and $B_{n}=\frac{B_{2n+1}^{\ast \ast
			}-B_{2n}^{\ast \ast }}{2}$ by Theorems \ref{teo2} and
		\ref{teo5}, we deduce that $B_{n}^{\ast }=\frac{3B_{2n+1}^{\ast
				\ast }-3B_{2n}^{\ast \ast }}{2}$. The others can be proved
		similarly.
	\end{proof}
	
	Thus we construct a one-to-one correspondence between all almost
	ba\-lan\-cing numbers of first type and all almost balancing numbers of
	second type.
	
	\section{Relationship with Pell Numbers.}
	
	In this section, we consider the relationship between all almost
	balancing numbers of first and second type and Pell numbers. It is known
	that the general terms of all balancing numbers can be given in terms of
	Pell numbers, namely 
	\begin{equation}
		B_{n}=\frac{P_{2n}}{2},b_{n}=\frac{P_{2n-1}-1}{2}%
		,C_{n}=P_{2n}+P_{2n-1},c_{n}=P_{2n-1}+P_{2n-2}  \label{pellk}
	\end{equation}
	for $n\geq 1$. Similarly we can give the following theorem.
	
	\begin{theorem}
		\label{teo7} The general terms of all almost balancing numbers of first type
		are 
		\begin{eqnarray*}
			&&B_{n}^{\ast }=\frac{3P_{2n}}{2},\ b_{2n}^{\ast }=\frac{4P_{2n}+P_{2n-1}-1}{%
				2},\ C_{n}^{\ast }=3P_{2n}+3P_{2n-1}, \\
			&&c_{2n-1}^{\ast }=5P_{2n-1}+P_{2n-2},\ c_{2n}^{\ast }=3P_{2n+1}-P_{2n}
		\end{eqnarray*}%
		for $n\geq 1,$ and $b_{2n-1}^{\ast }=\frac{8P_{2n-2}+3P_{2n-3}-1}{2}$ for $%
		n\geq 2,$ and the general terms of all almost balancing numbers of second
		type are%
		\begin{eqnarray*}
			&&B_{2n}^{\ast \ast }=\frac{P_{2n}+2P_{2n-1}}{2},\ b_{n}^{\ast \ast }=\frac{%
				3P_{2n-1}-1}{2} \\
			&&C_{2n}^{\ast \ast }=3P_{2n}-P_{2n-1},\ c_{n}^{\ast \ast
			}=3P_{2n-1}+3P_{2n-2}
		\end{eqnarray*}%
		for $n\geq 1,$ and $B_{2n-1}^{\ast \ast }=\frac{3P_{2n-2}+2P_{2n-3}}{2}%
		,C_{2n-1}^{\ast \ast }=5P_{2n-2}+P_{2n-3}$ for $n\geq 2$.
	\end{theorem}
	
	\begin{proof}
		Note that $B_{n}^{\ast }=3B_{n}$ and $B_{n}=\frac{P_{2n}}{2}$. So $%
		B_{n}^{\ast }=\frac{3P_{2n}}{2}$. Since $B_{2n-1}^{\ast \ast
		}=B_{n-1}$ $+C_{n-1}$ by Theorem \ref{teo2} and
		$B_{n}=\frac{P_{2n}}{2},C_{n}=P_{2n}+P_{2n-1}$ by (\ref{pellk}), we easily get
		\begin{equation*}
			B_{2n-1}^{\ast \ast }=\frac{P_{2n-2}}{2}+P_{2n-2}+P_{2n-3}=\frac{%
				3P_{2n-2}+2P_{2n-3}}{2}
		\end{equation*}%
		for $n\geq 2$ as we wanted. The other cases can be proved
		similarly.
	\end{proof}
	
	In Theorem \ref{teo7}, we can give the general terms of all almost balancing
	numbers of first and second type in terms of Pell numbers. Conversely, we
	can give the general terms of Pell numbers in terms of almost balancing
	numbers of first and second type as follows:
	
	\begin{theorem}
		\label{teo8} The general terms of Pell numbers are $P_{2n}=\frac{2B_{n}^{\ast }}{3}$ and $P_{2n-1}=b_{2n-1}^{\ast
			}-b_{2n-2}^{\ast }$ for $n\geq 1,$ or $P_{2n}=B_{2n+1}^{\ast \ast }-B_{2n}^{\ast \ast }$ and $P_{2n-1}=\frac{2b_{n}^{\ast \ast }+1}{3}$ 	for $n\geq 1$.
	\end{theorem}
	
	\begin{proof}
		It can be easily deduced from Theorem \ref{teo7}.
	\end{proof}
	
	Thus we construct a one-to-one correspondence between all almost
	ba\-lan\-cing numbers of first and second type and Pell numbers.

	\section{Concluding Remarks}
	
	For almost balancing and almost Lucas-balancing numbers of first and second
	type, in \cite{almost} Panda and Panda proved in Theorem 3.1 that the
	solutions of the Diophantine equation $8x^{2}+9=y^{2}$ in positive integers
	are given by $x=3B_{n}$ and $y=3C_{n}$ for $n\geq 1$. Similarly they proved
	in Theorem 3.2 that the solutions of the Diophantine equation $%
	8x^{2}-7=y^{2} $ in positive integers constitute two classes: the first
	class is $(x,y)=(B_{n}-2B_{n-1},$ $C_{n}-2C_{n-1})$, and the
	second class is $(x,y)=(2B_{n}-B_{n-1},2C_{n}-C_{n-1})$, for $n\geq 1$. Since 
	$B_{n}=\frac{\alpha ^{2n}-\beta ^{2n}}{4\sqrt{2}}$ and $C_{n}=\frac{\alpha
		^{2n}+\beta ^{2n}}{2}$, we easily deduce that 
	\begin{align*}
		B_{n}-2B_{n-1}& =\frac{\alpha ^{2n}-\beta ^{2n}}{4\sqrt{2}}-2(\frac{\alpha
			^{2n-2}-\beta ^{2n-2}}{4\sqrt{2}}) \\
		& =\frac{\alpha ^{2n-2}(1+2\sqrt{2})+\beta ^{2n-2}(-1+2\sqrt{2})}{4\sqrt{2}}
		\\
		& =\frac{\alpha ^{2n-2}-\beta ^{2n-2}}{4\sqrt{2}}+\frac{\alpha ^{2n-2}+\beta
			^{2n-2}}{2} \\
		& =B_{n-1}+C_{n-1} \\
		& =B_{2n-1}^{\ast \ast }
	\end{align*}%
	for $n\geq 1$. Similarly it can be shown that:
	\begin{gather*}
		C_{n}-2C_{n-1}=8B_{n-1}+C_{n-1}=C_{2n-1}^{\ast \ast }, \qquad
		2B_{n}-B_{n-1}=-B_{n}+C_{n}=B_{2n}^{\ast \ast }, \\
		2C_{n}-C_{n-1}=8B_{n}-C_{n}=C_{2n}^{\ast \ast },
	\end{gather*}
	 for $n\geq 1$, that is, we
	get same result in Theorem \ref{teo4}. Similarly for the almost cobalancing
	numbers of first and second type in \cite{almosttez}, Panda proved in
	Theorem 4.3.1 that the values of $x$ satisfying the Diophantine equation $%
	8x^{2}+8x+9=y^{2}$ in positive integers partition in two classes. The first
	class is given by $U_{n}=$ $\frac{3B_{n}+B_{n-1}-1}{2}$ and the second class
	is $V_{n}=\frac{3B_{n}+B_{n+1}-1}{2}$ for $n\geq 1$. Here we notice that 
	\begin{align*}
		U_{n}& =\frac{3B_{n}+B_{n-1}-1}{2} \\
		& =\frac{3(\frac{\alpha ^{2n}-\beta ^{2n}}{4\sqrt{2}})+\frac{\alpha
				^{2n-2}-\beta ^{2n-2}}{4\sqrt{2}}-1}{2} \\
		& =\frac{6(\frac{\alpha ^{2n}-\beta ^{2n}}{4\sqrt{2}})-\frac{\alpha
				^{2n}+\beta ^{2n}}{2}-1}{2} \\
		& =\frac{6B_{n}-C_{n}-1}{2} \\
		& =4b_{n}-b_{n-1}+1 \\
		& =b_{2n-1}^{\ast }
	\end{align*}%
	and similarly it can be shown that $V_{n}=b_{2n}^{\ast }$ as we proved in
	Theorem \ref{teo4}. But he did not determine the general terms of almost
	Lucas-cobalancing numbers of first and second type. Apart from these in this paper,
	
	\begin{enumerate}
		\item we determined the general terms of almost Lucas-cobalancing numbers
		of first and second type in Theorem \ref{teo4}.
		
		\item we can give the general terms of all balancing numbers in terms of all
		almost balancing numbers of first and second type in Theorem \ref{teo5}.
		Thus we construct a one-to-one correspondence between all balancing
		numbers and all almost balancing numbers of first and second type.
		
		\item We can give the general terms of all almost balancing numbers of first
		type in terms of all almost balancing numbers of second type and conversely give the
		general terms of all almost balancing numbers of second type in terms of all
		almost balancing numbers of first type in Theorem \ref{teo6}. Thus, we
		construct a one-to-one correspondence between all almost balancing numbers
		of first type and of second type.
		
		\item We can give the general terms of all almost balancing numbers of first
		and second type in terms of Pell numbers in Theorem \ref{teo7} and conversely give the general terms of Pell numbers in terms of almost balancing numbers
		of first and second type in Theorem \ref{teo8}. Thus, we construct a
		one-to-one correspondence between all almost balancing numbers of first
		and second type and Pell numbers.
	\end{enumerate}


\EditInfo{February 25, 2019}{March 16, 2021}{Attila Bérczes}

\end{paper}